\documentclass{amsart}
\usepackage{amsmath, amssymb}
\newcommand{\C}{\mathbb{C}}
\newcommand{\D}{\mathbb{D}}
\newcommand{\R}{\mathbb{R}}
\newcommand{\cS}{\mathcal{S}}
\newcommand{\cB}{\mathcal{B}}
\DeclareMathOperator{\dens}{dens}
\DeclareMathOperator{\Sing}{Sing}
\newtheorem{thm}{Theorem}
\newtheorem{cor}[thm]{Corollary}

\theoremstyle{remark}
\newtheorem*{rem}{Remark}
\newtheorem*{ack}{Acknowledgments}
\begin{document}
\bibliographystyle{amsalpha}
\title[On a conditional theorem of Littlewood]{On the exceptional set in a conditional theorem of Littlewood}
\begin{abstract}
  In 1952, Littlewood stated a conjecture about the average growth of
  spherical derivatives of polynomials, and showed that it would imply
  that for entire function of finite order, ``most'' preimages of
  almost all points are concentrated in a small subset of the
  plane. In 1988, Lewis and Wu proved Littlewood's conjecture. Using
  techniques from complex dynamics, we construct entire functions of
  finite order with a bounded set of singular values for which the set
  of exceptional preimages is infinite, with logarithmically growing
  cardinality.
\end{abstract}
\author{Lukas Geyer}
\address{Montana State University \\ Department of Mathematical
  Sciences \\ Bozeman, MT 59717 \\ USA}
\email{geyer@math.montana.edu}
\date{}
\maketitle
\section{Introduction and Main Result}
For a meromorphic function $f$ let $f^\# (z) =
\frac{2|f'(z)|}{1+|f(z)|^2}$ denote the spherical derivative of
$f$, and let $\Sing(f^{-1})$ denote the set of singular values,
i.e.~the set of all critical and asymptotic values of $f$. Let $\cS$
be the class of entire functions with finitely many singular values,
and $\cB$ the class of entire functions with a bounded set of finite
singular values. For Borel sets $A\subseteq \C$ we write $\#A$ for the
cardinality of a set $A$, and $|A|$ for its (two-dimensional) Lebesgue
measure. If $A,B \subset \C$ are Borel sets with $|B| \in (0,\infty)$ we
define $\dens(A,B) = \frac{|A\cap B|}{|B|}$ as the density of $A$ in
$B$. We denote the unit disk by $\D$, and the disk of radius $r$ centered
at 0 by $\D_r$.

If $R$ is a rational function of degree $n$, an application of
the Cauchy-Schwarz inequality yields
\begin{equation}
\begin{split}
\iint_\D R^\#(z) dx \, dy &\le \left( \iint_\D R^\#(z)^2 dx \, dy \right)^{1/2}
\left( \iint_\D dx \, dy \right )^{1/2} \\ &\le (4\pi n)^{1/2} \pi^{1/2} =
2\pi \sqrt{n},
\end{split}
\end{equation}
since $R$ covers the sphere $n$ times, and the area of the sphere is
$4\pi$. For general rational functions this is asymptotically best
possible, but Littlewood conjectured that for polynomials this
estimate could be improved. More precisely, he conjectured the
following, which was later proved by Lewis and Wu, building on earlier
partial results by Eremenko and Sodin in \cite{EremenkoSodin1986}
and \cite{EremenkoSodin1987}.

\begin{thm}[Lewis, Wu \cite{LewisWu1988}]
\label{LittlewoodConj}
There exist
absolute constants $C$ and $\alpha > 0$ with
\begin{equation}
\iint_\D P^\#(z) dx \, dy \le C n^{1/2-\alpha}
\end{equation}
for all $n$ and all polynomials $P$ of degree $n$.
\end{thm}

In fact, Lewis and Wu originally showed that one can choose $\alpha =
2^{-264}$, and Eremenko proved in \cite{Eremenko1991} that one cannot
choose $\alpha$ arbitrarily close to $1/2$. In \cite{Eremenko2002} a
connection between the supremum of the values of $\alpha$ for which
the theorem holds and the universal integral means spectrum was
conjectured. This connection is stated as a corollary of the main
result in \cite{BeliaevSmirnov2005}, but the proof of
the corollary crucially uses an unpublished result by Binder and Jones
about ``strong fractal approximation''. Assuming this result, both
upper and lower estimates on the supremum of $\alpha$-values in the
theorem can be improved considerably, see
\cite{HedenmalmShimorin2005}, \cite{Beliaev2008}, and
\cite{BeliaevSmirnov2010}.

Littlewood showed that his conjecture would have a
curious implication for the value distribution of entire functions of
finite order. Roughly speaking, most values are taken in a very small
subset of the plane. Since the conjecture is now a theorem,
Littlewood's conditional theorem becomes a corollary.

\begin{cor}[Littlewood\cite{Littlewood1952}]
\label{LittlewoodCor}
Let $f$ be an entire function of finite order $\rho \in (0,\infty)$,
and let $\beta \in (0,\alpha)$, where $\alpha$ is the constant of
Theorem \ref{LittlewoodConj}. Then there exists a constant $C_1$ and
an open set $S\subset\C$ with 
$\dens(S, \D_{r}) \le C_1 r^{-2\rho \beta}$ for all $r>0$, such
that for almost all $w\in \C$ and all
$\epsilon>0$, there exists a constant $C_2$ such that the set 
$E_w = f^{-1}(w)\setminus S$ satisfies 
$\# (E_w \cap \D_r) \le C_2 r^{\rho - (\alpha -\beta)\rho + \epsilon}$
for all $r>1$.
\end{cor}

We call $E_w$ the set of \emph{exceptional preimages}. Since we expect
to have roughly $r^{\rho}$ preimages in $\D_r$ for a function of order
$\rho$ and a typical point $w$, the estimate on the cardinality of
$E_w \cap \D_r$ shows that most preimages of typical points lie in
$S$, whose Lebesgue density in $\D_r$ is decreasing with a power of the
radius $r$.  Obviously, meromorphic functions of finite order do not
have this property, as shown by the Weierstrass $\wp$-function.

The question how large the set of exceptional preimages can be is
related to Epstein's ``order conjecture'', the question whether the
order of an entire function $f\in \cS$ is invariant under topological
equivalence in the sense of Eremenko and Lyubich
\cite{EremenkoLyubich1992}. If the number of exceptional preimages $\#
E_w$ is uniformly bounded on every compact set $K \subset \C \setminus
\Sing(f^{-1})$, then $f$ has the ``area property'', i.e.,
$\iint_{f^{-1}(K)} \frac{dx \, dy}{1+|z|^2} < \infty$ for every such
compact set $K$. This in turn implies invariance of the order of $f$
under topological equivalence in the class $\cS$.  For more
background, technical details and a similar construction to the one in
this paper of a function $f\in\cB$ which does not have the area
property see Epstein and Rempe \cite{EpsteinRempe2013}. The order
conjecture has recently been disproved by Bishop \cite{Bishop2013}.

In this note we show that in the class $\cB$ the exceptional set can
indeed contain $\ge c \log r$ points, as made precise in the following
theorem.

\begin{thm}
\label{CentralThm}
For almost every $\rho \in \left({\log 2}/{\log 3},\infty\right)$
there exists a function $f\in\cB$ of order $\rho$ and a set $W$ of
positive measure such that for any constants $C, \delta>0$, any
Borel set $S \subset \C$ satisfying $\dens(S,\D_r) \le C
r^{-\delta}$ for every $r>0$, and every $w\in W$, the set
$E_w=f^{-1}(w)\setminus S$ satisfies
\begin{equation}
\label{ExceptionalLog}
\liminf_{r\to\infty} \frac{\# (E_w \cap \D_r)}{\log r} \ge
\frac\rho{\log 2}.
\end{equation}
Furthermore, for every $\epsilon > 0$ there exists a function $f \in
\cB$ of order $\rho \in (1/2, 1/2+\epsilon)$ and a set $W$ of
positive measure satisfying
(\ref{ExceptionalLog}) for every $w\in W$ under the same assumptions.
\end{thm}

\begin{rem}
Entire functions in $\cB$ always have order $\rho \ge 1/2$ (see
\cite{BergweilerEremenko1995} and \cite{Langley1995}).
It is possible that our construction may be tweaked to yield
the result for all $\rho > 1/2$, but it will never produce examples of
order $\rho = 1/2$.
\end{rem}

\begin{ack}
I would like to thank Dmitri Beliaev, Alexandre Eremenko, Lasse Rempe,
and the referee for valuable remarks and suggestions.
\end{ack}

\section{Proof}
The entire functions we use will be Poincar\'e functions of quadratic
polynomials at repelling fixed points. We obtain the exceptional
preimages as preimages of a rotation domain under the Poincar\'e
function. In order to get almost all $\rho > {\log 2}/{\log 3}$,
we use explicit polynomials with fixed Siegel disks; in order to
obtain $\rho$ arbitrarily close to $1/2$, we choose perturbations of
the Chebyshev polynomial $T(z)=z^2-2$ with periodic Siegel disks. For
background on complex dynamics see \cite{CarlesonGamelin1993}.

Let $P(w) = \lambda w + w^2$ with $\lambda = e^{2\pi i \gamma}$.  By a
classical result of Siegel \cite{Siegel1942}, for almost all $\gamma \in \R$
the function $P$ can be linearized near $0$, i.e.~there exists an
analytic linearizing map $h(z) = z+O(z^2)$ near 0 such that $P(h(z)) =
h(\lambda z)$. The power series of $h$ has a finite radius of
convergence $R>0$, and $h$ maps $\D_R$ conformally onto the \emph{Siegel
  disk} $V$ of $P$ centered at 0. The polynomial $P$ has another
finite fixed point at $z_0 = 1-\lambda$ with multiplier $\mu :=
P'(z_0) = 2-\lambda$. Since $|\mu|>1$, there exists a local
linearizing function $f(z) = z_0 + z + O(z^2)$ with $P(f(z)) = f(\mu
z)$. In this case the functional equation allows to extend $f$ to an
entire function of order $\rho = {\log 2}/{\log |\mu|}$
(see \cite[\S 48]{Valiron1913}), the \emph{Poincar\'e function} of $P$ at
$z_0$. We now fix such a function $f$ associated to a polynomial $P$
with a Siegel disk $V=h(\D_R)$ centered at $0$, as well as the sub-Siegel disk
$W := h(\D_{R/2})$.

Now let $C,\delta >0$ be constants, and let $S \subset \C$ be a Borel
set satisfying $\dens(S,\D_r) \le C r^{-\delta}$ for every $r>0$. In
the following we use $C_k$ for constants depending only on $f$ and $S$.

We will show that $f\in\cB$ and that the exceptional set $E_w =
f^{-1}(w) \setminus S$ satisfies the asymptotic estimate
(\ref{ExceptionalLog}) for almost every $w \in W$. Since $|\mu| =
|2-\lambda|$ attains almost every value in the interval $(1,3)$, this
proves the theorem.

The set of singular values $\Sing(f^{-1})$ equals the post-critical
set of the polynomial $P$, i.e., the closure of the forward orbit of
the critical point \cite[Proposition
4.2]{MihaljevicBrandtPeter2012}. Since the latter is contained in the
Julia set of $P$, it is a bounded set disjoint from the simply
connected Siegel disk $V$. This implies both that $f\in\cB$, and that
$f$ maps every component of $f^{-1}(V)$ conformally onto $V$.

The Koebe Distortion Theorem implies that there exists an absolute
constant $M$ such that
\begin{equation}
\label{Koebe}
\frac1M  \le \frac{\dens (g(A),g(W))}{\dens(A,W)} \le M 
\end{equation}
for all conformal maps $g:V \to \C$ and all Borel sets $A \subseteq
W$ of positive measure.

Let $U_0$ be a component of $f^{-1}(V)$, and let $U_k = \mu^k
U$. Then $f(U_k) = f(\mu^k U_0) = P^k (V) = V$, so $(U_k)$ is a
sequence of components of $f^{-1}(V)$. Let $W_k = f^{-1}(W) \cap U_k$
and $S_k = S \cap W_k$. Since $W_0$ is a Borel set of measure
$|W_0|>0$ contained in some disk $\D_{C_1}$, we get that $W_k = \mu^k
W_0$ is a Borel set of measure $|W_k| = \mu^{2k} |W_0|$ contained in
$\D_{\mu^{k} C_1}$, so
\[
\dens(S,W_k) = \dens(S_k, W_k) \le C_2 \dens(S_k, \D_{\mu^{k} C_1})
\le C_3 \mu^{-k \delta}
\]
for all $k$. Applying (\ref{Koebe}) to the
branch of $f^{-1}$ mapping $V$ to $U_k$ yields 
\[
\dens(f(S_k), W) \le
C_4 |\mu|^{-k \delta}.
\]
Setting 
\[
E :=
\bigcap_{n=1}^\infty \bigcup_{k=n} f(S_k),
\]
we get
\[
\dens\left(E, W\right) \le \sum_{k=n}^\infty C_4 |\mu|^{-k \delta}
\]
for every $n$, so $\dens(E,W)=0$, and hence $|E|=0$. This implies
that almost every $w\in W$ satisfies $g^{-1}(w) \cap S_k = \emptyset$
for all but finitely many indices $k$. Thus for almost every $w\in W$,
the exceptional set $E_w$ contains
$\mu^k z_w$ for some $z_w \ne 0$ and all $k\ge 0$, and hence
\begin{equation}
\liminf_{k\to\infty} \frac{\# (E_w \cap \D_r)}{\log r} \ge
\frac{1}{\log |\mu|} = \frac{\rho}{\log 2}.
\end{equation}

In order to produce examples of order arbitrarily close to
$1/2$, we need to modify the construction slightly. Instead of using
explicit polynomials with Siegel fixed points, we use perturbations of
the Chebyshev polynomial $T(z) = z^2 - 2$ which have cycles of Siegel
disks. Existence of these polynomials is well-known, but for the
convenience of the reader we give a sketch of the proof.

The intermediate value theorem shows that there exists a decreasing
sequence of real numbers $(a_n)$ with $a_n \to -2$ such that $Q_n(z) =
z^2 + a_n$ has a super-attracting periodic point, i.e., it satisfies
$Q_n^{q_n}(0)=0$ for some $q_n \ge 1$. Perturbing $Q_n$ and using the
implicit function theorem, we get a sequence of polynomials $R_n(z) =
z^2+b_n$ with $-2< b_n < a_n$ having a parabolic periodic point $z_n$
with multiplier $-1$, i.e., $R_n^{q_n}(z_n) = z_n$ and
$(R_n^{q_n})'(z_n) = -1$. (Essentially this is the well-known fact
that there are Feigenbaum bifurcations arbitrarily close to the
Chebyshev polynomial in the quadratic family.) By the same result of
Siegel that we used in the first part of the proof, there are numbers
$\gamma$ arbitrarily close to $1/2$ such that any analytic function
$F(z) = e^{2\pi i \gamma}z + O(z^2)$ is linearizable. Since the
multiplier is a non-constant analytic function of the parameter near
$b_n$, there exists $c_n \in \C$ with $|c_n - b_n| < \frac1n$ such
that $P_n(z)=z^2 + c_n$ has a periodic Siegel disk of period $q_n$. In
this way we have constructed a sequence of polynomials $P_n(z) = z^2 +
c_n$ with periodic Siegel disks and $c_n \to -2$.

The repelling fixed point $z=2$ of $T(z) = z^2-2$ varies analytically
with the parameter, so $P_n$ has a repelling fixed point $z_n$ with
$z_n \to 2$ and $\mu_n = P_n'(z_n) \to 4$ for $n\to\infty$. It
follows from classical results in complex dynamics that the Julia set
of $P_n$ is connected, and this implies $|\mu_n| <4$ (see \cite{Buff2003}).
Let $f_n$ denote the Poincar\'e function of $P_n$ at the fixed point
$z_n$. Since $P_n$ has connected Julia set, we get $f_n \in \cB$ with
order $\rho_n = {\log 2}/{\log \mu_n} \to 1/2$. It remains
to show that $f_n$ satisfies (\ref{ExceptionalLog}), and this follows
along the same lines as in the first part of the proof.

Let $n$ be fixed and, and let $U_1,\ldots U_{q-1}$ be
the cycle of Siegel disks of $P_n$ containing periodic points $\zeta_1,
\ldots, \zeta_{q-1}$. Let $h_k: \D_{R_k} \to U_k$ be
the linearizing map of $P^q$ in $U_k$, normalized as $h_k(z) = \zeta_k
+ z + O(z^2)$. Now we let $W = \bigcup_{k=1}^{q-1}
h_k(\D_{R_k/2})$. Then $P_n(W)=W$, and $W$ is a finite union of
sub-Siegel disks, so we also get (\ref{Koebe}) for all maps $g$
which are conformal in any $U_j$, where the constants do not depend on
$j$. Applying this to branches of $f_n^{-1}$ exactly as in the first
part of the proof yields the desired estimate (\ref{ExceptionalLog}).
\qed

\bibliography{MathSciNet}

\end{document}